# Some elementary aspects of q-Fibonacci and q-Lucas polynomials

Johann Cigler
johann.cigler@univie.ac.at

**Abstract.**
Based on well-known properties of Fibonacci and Lucas numbers and polynomials we give a self-contained approach to some bivariate $q$ – analogs.

## 0. Introduction

We give a simple approach to some elementary aspects of the $q$ – Fibonacci polynomials

$$F_n(x,s,q) = \sum_{k=0}^{\lfloor \frac{n-1}{2} \rfloor} q^{k^2} \begin{bmatrix} n-1-k \\ k \end{bmatrix}_q s^k x^{n-1-2k},$$ introduced by L. Carlitz [5], the polynomials

$$Fib_n(x,s,q) = \sum_{k=0}^{\lfloor \frac{n-1}{2} \rfloor} q^{\binom{k+1}{2}} \begin{bmatrix} n-1-k \\ k \end{bmatrix}_q s^k x^{n-1-2k},$$ which have been studied in [9], and of the

associated $q$ – Lucas polynomials.

We restrict ourselves to typical properties such as recurrence relations, closed formulas, combinatorial interpretations, and generating functions.

We emphasize the analogies with Fibonacci and Lucas numbers and polynomials by recalling some well-known properties of them and give a brief introduction to the required tools from $q$ – analysis.

There is a **notational dilemma** with regard to Fibonacci numbers.

The Fibonacci Quarterly, a journal devoted to Fibonacci numbers and their properties and generalizations, insists that the Fibonacci numbers $F_n$ are defined by $F_{n+2} = F_{n+1} + F_n$, $F_0 = 0$, $F_1 = 1$. This notation is adapted to **arithmetical properties** such as Binet's formula.

Some combinatorialists like Herbert S. Wilf [18], on the other hand, choose as initial values $F_0 = 1$ and $F_1 = 1$, which is more appropriate for **combinatorial interpretations.**

We will need both versions. In order to avoid misunderstandings, we shall use $f_n$ for the combinatorial version and $F_n$ for the usual one. Thus $(F_n)_{n \geq 0} = (0,1,1,2,3,5,8,\cdots)$ and $(f_n)_{n \geq 0} = (1,1,2,3,5,8,\cdots)$.

A somewhat different situation occurs for the Lucas numbers. For arithmetical properties the sequence $(L_n)_{n \geq 0} = (2,1,3,4,7,)$ is the right one. For combinatorial interpretations the sequence $(l_n)_{n \geq 0} = (1,1,3,4,7,)$ is better suited.

Similar notational distinctions will be made for Fibonacci and Lucas polynomials and their $q$ – analogs.

## 1. Fibonacci and Lucas numbers revisited

The **Fibonacci numbers** $F_n$ are recursively defined by $F_n = F_{n-1} + F_{n-2}$ with initial values $F_0 = 0$ and $F_1 = 1$. The first terms are $(F_n)_{n \geq 0} = (0, 1, 1, 2, 3, 5, 8, 13, 21, \cdots)$.

The **Lucas numbers** $L_n$ satisfy the same recurrence relation $L_n = L_{n-1} + L_{n-2}$ but with initial values $L_0 = 2$ and $L_1 = 1$. The first terms are $(L_n)_{n \geq 0} = (2, 1, 3, 4, 7, 11, 18, \cdots)$.

Let $\alpha = \dfrac{1+\sqrt{5}}{2}$ and $\beta = \dfrac{1-\sqrt{5}}{2}$ be the roots of $z^2 - z - 1 = 0$.

Then each linear combination $a\alpha^n + b\beta^n$ satisfies the same recurrence. Choosing $a, b$ such that the initial values are satisfied gives the **Binet formulas**

$$F_n = \frac{\alpha^n - \beta^n}{\alpha - \beta}, \qquad (1)$$
$$L_n = \alpha^n + \beta^n.$$

As a simple application let us show that $F_{kn}$ is a multiple of $F_n$ for each $k \in \mathbb{N}$.

Binet's formula gives $F_{kn} = \dfrac{\alpha^{kn} - \beta^{kn}}{\alpha - \beta} = \dfrac{\alpha^{kn} - \beta^{kn}}{\alpha^n - \beta^n} \cdot \dfrac{\alpha^n - \beta^n}{\alpha - \beta} = F_k^{(n)} F_n$.

The numbers $F_k^{(n)} = \dfrac{\alpha^{kn} - \beta^{kn}}{\alpha^n - \beta^n}$ satisfy $F_k^{(n)} = L_n F_{k-1}^{(n)} + F_{k-2}^{(n)}$ with $F_0^{(n)} = 0$ and $F_1^{(n)} = 1$ because $(z - \alpha^n)(z - \beta^n) = z^2 - L_n z - 1$.

The Binet formulas can be used to extend $F_n$ and $L_n$ to negative $n$. Observing that $\alpha\beta = -1$ this gives $F_{-n} = (-1)^{n-1} F_n$ and $L_{-n} = (-1)^n L_n$. We will restrict ourselves to $n \geq 0$, but shall use for some formulas the fact that $F_{-1} = 1$.

For example, it is convenient to have the well-known formula
$$L_n = F_{n+1} + F_{n-1} \qquad (2)$$
for all $n \in \mathbb{N}$.

It follows from
$$\frac{\alpha^{n+1} - \beta^{n+1}}{\alpha - \beta} + \frac{\alpha^{n-1} - \beta^{n-1}}{\alpha - \beta} = \frac{\alpha^{n+1} - \beta^{n+1}}{\alpha - \beta} - \alpha\beta \frac{\alpha^{n-1} - \beta^{n-1}}{\alpha - \beta} = \alpha^n + \beta^n.$$

The recurrence for the Fibonacci numbers can be written in matrix form as

$$\begin{pmatrix} F_{n-1} & F_n \\ F_n & F_{n+1} \end{pmatrix} = \begin{pmatrix} 0 & 1 \\ 1 & 1 \end{pmatrix} \begin{pmatrix} F_{n-2} & F_{n-1} \\ F_{n-1} & F_n \end{pmatrix} = \begin{pmatrix} 0 & 1 \\ 1 & 1 \end{pmatrix}^n = U^n \qquad (3)$$

with $U = \begin{pmatrix} 0 & 1 \\ 1 & 1 \end{pmatrix} = \begin{pmatrix} F_0 & F_1 \\ F_1 & F_2 \end{pmatrix}$.

From $U^{m+n} = U^m U^n$ we get $\begin{pmatrix} F_{m+n-1} & F_{m+n} \\ F_{m+n} & F_{m+n+1} \end{pmatrix} = \begin{pmatrix} F_{m-1} & F_m \\ F_m & F_{m+1} \end{pmatrix} \begin{pmatrix} F_{n-1} & F_n \\ F_n & F_{n+1} \end{pmatrix}$

and thus
$$F_{m+n} = F_{m-1} F_n + F_m F_{n+1}. \tag{4}$$

Taking determinants, we get **Cassini's identity**
$$F_{n-1} F_{n+1} - F_n^2 = (-1)^n. \tag{5}$$

The identity
$$U^2 = \begin{pmatrix} 1 & 1 \\ 1 & 2 \end{pmatrix} = \begin{pmatrix} 1 & 0 \\ 0 & 1 \end{pmatrix} + \begin{pmatrix} 0 & 1 \\ 1 & 1 \end{pmatrix} = I + U$$
implies
$$U^{2n+m} = (I+U)^n U^m = \sum_{k=0}^n \binom{n}{k} U^{k+m}$$
which gives the **doubling formula**
$$F_{2n+m} = \sum_{k=0}^n \binom{n}{k} F_{k+m}, \tag{6}$$

Let us also note that the Lucas numbers $L_n$ are the trace of $U^n$
$$L_n = F_{n+1} + F_{n-1} = tr U^n. \tag{7}$$
Since
$$(1-z-z^2) \sum_{n \geq 0} F_{n+1} z^n = \sum_{n \geq 0} F_{n+1} z^n - \sum_{n \geq 1} F_n z^n - \sum_{n \geq 2} F_{n-1} z^n = 1 + z - z + \sum_{n \geq 2} (F_{n+1} - F_n - F_{n-1}) z^n = 1$$
the generating function of the Fibonacci numbers is
$$\sum_{n \geq 0} F_{n+1} z^n = \frac{1}{1-z-z^2}. \tag{8}$$
For the Lucas numbers we get in the same way
$$\sum_{n \geq 0} L_n z^n = \frac{2-z}{1-z-z^2}. \tag{9}$$
From the generating function (8) we get
$$\sum_{n \geq 0} F_{n+1} z^n = \frac{1}{1-z(1+z)} = \sum_{\ell \geq 0} z^\ell (1+z)^\ell = \sum_{\ell,k \geq 0} \binom{\ell}{k} z^{\ell+k} = \sum_{n \geq 0} z^n \sum_{2k \leq n} \binom{n-k}{k}$$
which gives
$$F_{n+1} = \sum_{j=0}^{\lfloor \frac{n}{2} \rfloor} \binom{n-j}{j}. \tag{10}$$

A nice **combinatorial model for Fibonacci numbers** are sets of **Morse sequences** $c_1 c_2 \cdots c_k$ where the elements $c_i$ are dots $a := \bullet$ or dashes $b := -$.

We associate to a Morse sequence $c = c_1 c_2 \cdots c_k$ the length

$$L(c) = L(c_1 c_2 \cdots c_k) = \sum_{i=1}^{n} L(c_i) \tag{11}$$

with $L(a) = 1$ and $L(b) = 2$.

Let $M_n$ be the set of all Morse sequences $c_1 c_2 \cdots c_k$ of length $n$.

$M_1 = \{\bullet\}$, $M_2 = \{\bullet\bullet, -\}$, $M_3 = \{\bullet\bullet\bullet, \bullet-, -\bullet\}$, $M_4 = \{\bullet\bullet\bullet\bullet, \bullet\bullet-, \bullet-\bullet, -\bullet\bullet, --\}$.

The elements of $M_n$ can also be interpreted as coverings of the interval $[0,n]$ of the real line with intervals of length 1 (the dots) and intervals of length 2 (the dashes).

The number $|M_n|$ of elements of $M_n$ is

$$|M_n| = f_n. \tag{12}$$

For this is true for $n=1$ and $n=2$ because $|M_1| = 1 = f_1$ and $|M_2| = 2 = f_2$.

In the general case $M_n = M_{n,a} \cup M_{n,b}$, where $M_{n,a}$ consists of all sequences with $c_1 = a$ and $M_{n,b}$ consists of all sequences with $c_1 = b$. Obviously $|M_{n,a}| = |M_{n-1}|$ and $|M_{n,b}| = |M_{n-2}|$. Therefore we get $|M_n| = |M_{n-1}| + |M_{n-2}| = f_{n-1} + f_{n-2} = f_n$.

For $n = 0$ we let $M_0 = \{\varnothing\}$ be the set consisting of the empty sequence which gives $|M_0| = 1 = f_0$.

From a combinatorial point of view formula (10) states that for each $j$ there are $\binom{n-j}{j}$ Morse sequences $c = c_1 c_2 \cdots c_{n-j}$ of length $n$ whose $j$ dashes can be distributed in $\binom{n-j}{j}$ ways.

The (combinatorial) **Lucas numbers** $l_n$ can be interpreted as weights of **periodic coverings of the real line with dots and dashes.**

Let $M_n^*$ be the set of coverings with period $n$. These are uniquely determined by their restrictions to $[0, n]$. Those coverings where $[0,1]$ is covered by a dot or the first part of a dash can be identified with $M_n$. Those where $[0,1]$ is covered by the second part of a dash can be identified with all linear coverings of $[1, n-1]$, i.e. with $M_{n-2}$.

Therefore, we get for $n \geq 2$

$$|M_n^*| = |M_n| + |M_{n-2}| = f_n + f_{n-2} = l_n. \tag{13}$$

It is clear that $|M_1^*| = 1 = l_1$. For $n = 0$ we get the empty covering which gives $|M_0^*| = 1 = l_0$.

Another formula for the Lucas numbers for $n > 0$ is

$$L_n = l_n = \sum_j \binom{n-j}{j} \frac{n}{n-j}. \tag{14}$$

This follows from

$$L_n = f_n + f_{n-2} = \sum_j \left(\binom{n-j}{j} + \binom{n-2-j}{j}\right) = \sum_j \left(\binom{n-j}{j} + \binom{n-1-j}{j-1}\right)$$

$$= \sum_j \binom{n-j}{j}\left(1 + \frac{j}{n-j}\right) = \sum_j \binom{n-j}{j}\frac{n}{n-j}.$$

## 2. Bivariate Fibonacci and Lucas polynomials

Let us recall some well-known results about Fibonacci and Lucas polynomials.

By comparing coefficients it is easily verified that the **(arithmetical) Fibonacci polynomials**

$$F_n(x,s) = \sum_{j=0}^{\lfloor \frac{n-1}{2} \rfloor} \binom{n-1-j}{j} s^j x^{n-1-2j} \tag{15}$$

satisfy

$$F_n(x,s) = xF_{n-1}(x,s) + sF_{n-2}(x,s) \tag{16}$$

with initial values $F_0(x,s) = 0$ and $F_1(x,s) = 1$.

The corresponding Binet formula is

$$F_n(x,s) = \frac{\alpha(x,s)^n - \beta(x,s)^n}{\alpha(x,s) - \beta(x,s)} \tag{17}$$

where

$$\alpha(x,s) = \frac{x + \sqrt{x^2 + 4s}}{2},$$
$$\beta(x,s) = \frac{x - \sqrt{x^2 + 4s}}{2} \tag{18}$$

are the roots of $z^2 - xz - s = 0$.

If we define the weight of a Morse sequence $c = c_1 c_2 \cdots c_k$ by $w(c_1 c_2 \cdots c_k) = \prod_{j=1}^{k} w(c_i)$ with $w(a) = x$ and $w(b) = s$ then the **bivariate (combinatorial) Fibonacci polynomials**

$$f_n(x,s) = \sum_{j=0}^{\lfloor \frac{n}{2} \rfloor} \binom{n-j}{j} s^j x^{n-2j} \tag{19}$$

can be interpreted as the weight of the set $M_n$ of all Morse sequences of length $n$

$$f_n(x,s) = w(M_n) = \sum_{c \in M_n} w(c). \tag{20}$$

For there are $\binom{n-j}{j}$ Morse sequences $c_1 c_2 \cdots c_{n-j}$ of length $n$ with $j$ dashes and the weight of each sequence is $w(c_1 c_2 \cdots c_{n-j}) = s^j x^{n-2j}$.

Thus $f_n(x,s)$ is a polynomial in $x$ and $s$ with $\deg_x f_n(x,s) = n.$
The recurrence
$$f_n(x,s) = xf_{n-1}(x,s) + sf_{n-2}(x,s) \tag{21}$$
can also be verified by observing that $w(M_{n,a}) = xw(M_{n-1}) = xf_{n-1}(x,s)$ and $w(M_{n,b}) = sw(M_{n-2}) = sf_{n-2}(x,s).$

For arbitrary $x, y$ we get (cf. [13]
$$F_n(x+y, -xy) = \frac{x^n - y^n}{x - y}, \tag{22}$$
because $\dfrac{x^n - y^n}{x - y} = (x+y)\dfrac{x^{n-1} - y^{n-1}}{x - y} - xy\dfrac{x^{n-2} - y^{n-2}}{x - y}.$

Let us mention the special case $x = \dfrac{1 + \sqrt{-3}}{2} = e^{\frac{2\pi i}{6}}$ and $y = \dfrac{1 - \sqrt{-3}}{2} = e^{-\frac{2\pi i}{6}}$, which implies that $F_n(1,-1) = \sum_{j=0}^{\lfloor \frac{n-1}{2} \rfloor} (-1)^j \binom{n-1-j}{j}$ is periodic with period 6:

$$\begin{aligned} F_{3n}(1,-1) &= 0, \\ F_{3n+1}(1,-1) &= (-1)^n, \\ F_{3n+2}(1,-1) &= (-1)^n. \end{aligned} \tag{23}$$

The recurrence for the Fibonacci polynomials can be written in matrix form as
$$\begin{pmatrix} sF_{n-1}(x,s) & F_n(x,s) \\ sF_n(x,s) & F_{n+1}(x,s) \end{pmatrix} = \begin{pmatrix} 0 & 1 \\ s & x \end{pmatrix} \begin{pmatrix} sF_{n-2}(x,s) & F_{n-1}(x,s) \\ sF_{n-1}(x,s) & F_n(x,s) \end{pmatrix} = \begin{pmatrix} 0 & 1 \\ s & x \end{pmatrix}^n = U(x,s)^n \tag{24}$$
with $U(x,s) = \begin{pmatrix} 0 & 1 \\ s & x \end{pmatrix}.$

Taking determinants, we get **Cassini's identity**
$$F_n^2(x,s) - F_{n-1}(x,s)F_{n+1}(x,s) = (-s)^{n-1}. \tag{25}$$
From $U(x,s)^{m+n} = U(x,s)^m U(x,s)^n$ we get
$$F_{m+n}(x,s) = F_m(x,s)F_{n+1}(x,s) + sF_{m-1}(x,s)F_n(x,s). \tag{26}$$

$$U(x,s)^2 = \begin{pmatrix} s & x \\ sx & x^2+s \end{pmatrix} = s\begin{pmatrix} 1 & 0 \\ 0 & 1 \end{pmatrix} + x\begin{pmatrix} 0 & 1 \\ s & x \end{pmatrix} = sI + xU(x,s)$$
implies
$$U(x,s)^{2n+m} = (sI + xU(x,s))^n U(x,s)^m = \sum_{k=0}^n \binom{n}{k} s^k x^{n-k} U(x,s)^{n-k+m},$$

which gives the **doubling formula**

$$F_{2n+m}(x,s) = \sum_{k=0}^{n} \binom{n}{k} s^k x^{n-k} F_{n-k+m}(x,s). \tag{27}$$

The **bivariate Lucas polynomials** are defined by

$$L_n(x,s) = F_{n+1}(x,s) + sF_{n-1}(x,s) = \sum_{j=0}^{\lfloor \frac{n}{2} \rfloor} \binom{n-j}{j} \frac{n}{n-j} s^j x^{n-2j} \quad \text{for } n > 0, \tag{28}$$

$$L_0(x,s) = 2.$$

They satisfy

$$L_n(x,s) = xL_{n-1}(x,s) + sL_{n-2}(x,s) \tag{29}$$
$$\text{with } L_0(x,s) = 2 \text{ and } L_1(x,s) = x.$$

The Binet formula for the Lucas polynomials gives

$$L_n(x,s) = \alpha(x,s)^n + \beta^n(x,s), \tag{30}$$

**Proof.**

$$L_n(x,s) = F_{n+1}(x,s) + sF_{n-1}(x,s) = \frac{\alpha(x,s)^{n+1} - \beta(x,s)^{n+1}}{\alpha(x,s) - \beta(x,s)} - \alpha(x,s)\beta(x,s) \frac{\alpha(x,s)^{n-1} - \beta(x,s)^{n-1}}{\alpha(x,s) - \beta(x,s)}$$

$$= \frac{1}{\alpha(x,s) - \beta(x,s)} \left( \alpha(x,s)^{n+1} - \alpha(x,s)^n \beta(x,s) + \alpha(x,s)\beta^n(x,s) - \beta(x,s)^{n+1} \right)$$

$$= \alpha(x,s)^n + \beta^n(x,s)$$

Note that

$$L_n(x,s) = trU(x,s)^n. \tag{31}$$

We define the **bivariate (combinatorial) Lucas polynomials** $l_n(x,s)$ as the weight of $M_n^*$,

$$l_n(x,s) = w(M_n^*). \tag{32}$$

With the same argument as above, we get for $n \geq 2$

$$w(M_n^*) = w(M_n) + sw(M_{n-2}) = f_n(x,s) + sf_{n-2}(x,s). \tag{33}$$

Also note that $w(M_1^*) = x$. In order that each $l_n(x,s)$ is a monic polynomial of degree $n$ in $x$ we set $l_0(x,s) = 1$.

Thus we get

$$l_n(x,s) = f_n(x,s) + sf_{n-2}(x,s) \tag{34}$$

with initial values $l_0(x,s) = 1$ and $l_1(x,s) = x$.

The polynomials $l_n(x,y)$ are monic of degree $n$. They satisfy

$$l_n(x,s) = xl_{n-1}(x,s) + t_{n-2}(s)l_{n-2}(x,s) \tag{35}$$
$$\text{with } t_0(s) = 2s \text{ and } t_n(s) = s \text{ for } n > 0.$$

An analog of (22) is
$$L_n(x+y,-xy) = x^n + y^n. \tag{36}$$
This is true for $n=0$ and $n=1$. In the general case we have
$$x^{n+1} + y^{n+1} = (x+y)(x^n + y^n) - xy(x^{n-1} + y^{n-1}).$$
In the same way as above, we get the generating functions

$$\sum_{n\geq 0} f_n(x,s) z^n = \frac{1}{1-xz-sz^2} = \frac{1}{1-xz} \frac{1}{1 - \frac{sz^2}{1-xz}} = \sum_{n\geq 0} \frac{s^n}{(1-xz)^{n+1}} z^{2n},$$

$$\sum_{n\geq 0} L_n(x,s) z^n = \frac{2-xz}{1-xz-sz^2} = (2-xz) \sum_{n\geq 0} \frac{s^n}{(1-xz)^{n+1}} z^{2n}. \tag{37}$$

For the combinatorial Fibonacci and Lucas polynomials there exist **nice inversion formulas**.

$$\sum_{k=0}^{\lfloor n/2 \rfloor} (-s)^k \binom{n}{k} l_{n-2k}(x,s) = x^n,$$

$$\sum_{k=0}^{\lfloor n/2 \rfloor} (-s)^k \left( \binom{n}{k} - \binom{n}{k-1} \right) f_{n-2k}(x,s) = x^n. \tag{38}$$

**Proof**
These are essentially the Chebyshev inverse relations (cf. [14]). Let us give a direct proof. We consider the cases for odd and even $n$ separately. First we prove the formula for the Lucas polynomials. Noting that $\alpha(x,s)\beta(x,s) = -s$ we get

$$\sum_{k=0}^{n} (-s)^k \binom{2n+1}{k} l_{2n+1-2k}(x,s) = \sum_{k=0}^{n} (\alpha(x,s)\beta(x,s))^k \binom{2n+1}{k} \left( \alpha(x,s)^{2n+1-2k} + \beta(x,s)^{2n+1-2k} \right)$$

$$= \sum_{k=0}^{n} \binom{2n+1}{k} \left( \alpha(x,s)^{2n+1-k} \beta(x,s)^k + \alpha(x,s)^k \beta(x,s)^{2n+1-k} \right) = \sum_{k=0}^{2n+1} \binom{2n+1}{k} \alpha(x,s)^k \beta(x,s)^{2n+1-k}$$

$$= (\alpha(x,s) + \beta(x,s))^{2n+1} = x^{2n+1}$$

and

$$\sum_{k=0}^{n} (-s)^k \binom{2n}{k} l_{2n-2k}(x,s) = \sum_{k=0}^{n-1} (\alpha(x,s)\beta(x,s))^k \binom{2n}{k} \left( \alpha(x,s)^{2n-2k} + \beta(x,s)^{2n-2k} \right) + (-s)^n \binom{2n}{n}$$

$$= \sum_{k=0}^{n-1} \binom{2n}{k} \left( \alpha(x,s)^{2n-k} \beta(x,s)^k + \alpha(x,s)^k \beta(x,s)^{2n-k} \right) + (\alpha(x,s)\beta(x,s))^n \binom{2n}{n} = \sum_{k=0}^{2n} \binom{2n}{k} \alpha(x,s)^k \beta(x,s)^{2n-k}$$

$$= (\alpha(x,s) + \beta(x,s))^{2n} = x^{2n}.$$

The result for the Fibonacci polynomials follows from

$$x^n = \sum_{k=0}^{\lfloor n/2 \rfloor} (-s)^k \binom{n}{k} l_{n-2k}(x,s) = \sum_{k=0}^{\lfloor n/2 \rfloor} (-s)^k \binom{n}{k} \left(f_{n-2k}(x,s) + sf_{n-2k-2}(x,s)\right)$$

$$= \sum_{k=0}^{\lfloor n/2 \rfloor} (-s)^k \binom{n}{k} f_{n-2k}(x,s) - \sum_{k=0}^{\lfloor n/2 \rfloor} (-s)^{k+1} \binom{n}{k} f_{n-2k-2}(x,s)$$

$$= \sum_{k=0}^{\lfloor n/2 \rfloor} (-s)^k f_{n-2k}(x,s) \left(\binom{n}{k} - \binom{n}{k-1}\right)$$

**Remark**

A sequence of monic polynomials $p_n(x)$ of degree $n$ is called orthogonal if

$$p_n(x) = (x - s_{n-1})p_{n-1}(x) - t_{n-2}p_{n-2}(x) \tag{39}$$

for some values $s_n$ and $t_n$. Let $\Lambda = \Lambda_{p(x)}$ denote the linear functional on the polynomials defined by $\Lambda(p_n(x)) = [n = 0]$. Then (39) implies the orthogonality relations $\Lambda(p_n(x)p_m(x)) = 0$ for $n \neq m$ and $\Lambda(p_n(x)^2) \neq 0$. The numbers $\Lambda(x^n)$ are called moments of $\Lambda$.

Equations (21) and (35) show that the polynomials $f_n(x,s)$ and $l_n(x,s)$ are orthogonal.

By (38) the moments of the Fibonacci and Lucas polynomials are

$$\Lambda_{f(x,s)}\left(x^{2n}\right) = (-s)^n C_n = (-s)^n \frac{1}{n+1}\binom{2n}{n}, \quad \Lambda_{f(x,s)}\left(x^{2n+1}\right) = 0,$$

$$\Lambda_{l(x,s)}\left(x^{2n}\right) = (-s)^n \binom{2n}{n}, \quad \Lambda_{l(x,s)}\left(x^{2n+1}\right) = 0. \tag{40}$$

## 2. Some definitions and results from q-analysis

### 2.1. Preliminaries

Let $q \neq 1$ be a real number and let $[n]_q = \dfrac{1-q^n}{1-q}$, $[n]_q! = [1]_q[2]_q \cdots [n]_q$ with $0_q! = 1$ and $\begin{bmatrix}n\\k\end{bmatrix}_q = \dfrac{[n]_q!}{[k]_q![n-k]_q!}$ for $0 \leq k \leq n$ and $\begin{bmatrix}n\\k\end{bmatrix}_q = 0$ else.

For $q \to 1$ we have $[n]_q \to n$ and $\begin{bmatrix}n\\k\end{bmatrix}_q \to \binom{n}{k}$. Therefore $[n]_q$ is called a $q$-analogue of $n$ and $\begin{bmatrix}n\\k\end{bmatrix}_q$ a $q$-analogue of $\binom{n}{k}$.

The $q$-binomial coefficients satisfy

$$\begin{bmatrix} n \\ k \end{bmatrix}_q = \begin{bmatrix} n \\ n-k \end{bmatrix}_q \quad (41)$$

and

$$\begin{bmatrix} n+1 \\ k \end{bmatrix}_q = q^k \begin{bmatrix} n \\ k \end{bmatrix}_q + \begin{bmatrix} n \\ k-1 \end{bmatrix}_q,$$

$$\begin{bmatrix} n+1 \\ k \end{bmatrix}_q = \begin{bmatrix} n \\ k \end{bmatrix}_q + q^{n+1-k} \begin{bmatrix} n \\ k-1 \end{bmatrix}_q. \quad (42)$$

Let $D_q$ be the $q$-differentiation operator on the vector space of polynomials in $x$ defined by $D_q p(x) = \dfrac{p(x) - p(qx)}{(1-q)x}$. It is uniquely determined by $D_q x^n = [n]_q x^{n-1}$ for all $n \in \mathbb{N}$. For $q \to 1$ it converges to ordinary differentiation.

### 2.2. Some q-analogs of the binomial theorem

**Binomial theorem for q-commuting operators ([6],[16])**
*Let $A$ and $B$ be linear operators on the polynomials satisfying*

$$BA = qAB \quad (43)$$

*then*

$$(A+B)^n = \sum_{k=0}^n \begin{bmatrix} n \\ k \end{bmatrix} A^k B^{n-k}. \quad (44)$$

**Proof**
By induction we get

$$(A+B)\sum_{k=0}^n \begin{bmatrix} n \\ k \end{bmatrix} A^k B^{n-k} = \sum_{k=0}^n \begin{bmatrix} n \\ k \end{bmatrix} A^{k+1} B^{n-k} + \sum_{k=0}^n \begin{bmatrix} n \\ k \end{bmatrix} q^k A^k B^{n-k+1}$$

$$= \sum_{k=0}^{n+1} \begin{bmatrix} n \\ k-1 \end{bmatrix} A^k B^{n+1-k} + \sum_{k=0}^{n+1} \begin{bmatrix} n \\ k \end{bmatrix} q^k A^k B^{n-k+1} = \sum_{k=0}^{n+1} \left( \begin{bmatrix} n \\ k-1 \end{bmatrix} + q^k \begin{bmatrix} n \\ k \end{bmatrix} \right) A^k B^{n-k+1} = \sum_{k=0}^{n+1} \begin{bmatrix} n+1 \\ k \end{bmatrix} A^k B^{n+1-k}.$$

Let $\varepsilon = \varepsilon_q$ be the linear operator on the polynomials defined by $\varepsilon_q p(x) = p(qx)$ and let $M$ be the multiplication operator with $x$. Then $\varepsilon_q M = qM\varepsilon_q$ because $\varepsilon_q M p(x) = \varepsilon_q x p(x) = qx p(qx) = qM\varepsilon_q p(x)$. Since there can be no misunderstandings in the following we simply write $x$ in place of $M$.

If we apply the operator $(x + y\varepsilon_q)^n$ to the constant polynomial $p(x) = 1$ we get the **Rogers-Szegö polynomials**

$$R_n(x,y,q) = (x + y\varepsilon_q)^n 1 = \sum_{k=0}^n \begin{bmatrix} n \\ k \end{bmatrix}_q x^k y^{n-k}. \quad (45)$$

They satisfy

$$D_q R_n(x,y,q) = [n]_q R_{n-1}(x,y,q) \tag{46}$$

and

$$R_n(x,y,q) = (x+y) R_{n-1}(x,y,q) + (q^{n-1}-1) xy R_{n-2}(x,y,q). \tag{47}$$

For

$$D_q R_n(x,y,q) = \sum_{k=0}^{n} \begin{bmatrix} n \\ k \end{bmatrix}_q [k]_q x^{k-1} y^{n-k} = [n]_q \sum_{k=0}^{n} \begin{bmatrix} n-1 \\ k-1 \end{bmatrix}_q x^{k-1} y^{n-k} = [n]_q \sum_{k=0}^{n-1} \begin{bmatrix} n-1 \\ k \end{bmatrix}_q x^k y^{n-1-k} = [n]_q R_{n-1}(x,y,q)$$

and

$$\varepsilon_q = 1 + (q-1) x D_q. \tag{48}$$

The last identity holds because
$$(1+(q-1)xD_q) x^n = x^n + (q-1)x[n]x^{n-1} = x^n (1+q^n-1) = (qx)^n = \varepsilon_q x^n.$$

This finally gives the recurrence

$$R_n(x,y,q) = (x + y(1+(q-1)xD_q)) R_{n-1}(x,y,q) = (x+y) R_{n-1}(x,y,q) + (q^{n-1}-1) xy R_{n-2}(x,y,q).$$

By setting $A = x\varepsilon_q$ and $B = y\varepsilon_q$ we get Rothe's theorem

$$(y+x)(y+qx)\cdots(y+q^{n-1}x) = \sum_{k=0}^{n} q^{\binom{k}{2}} \begin{bmatrix} n \\ k \end{bmatrix} x^k y^{n-k}. \tag{49}$$

For $\varepsilon(x\varepsilon) = q(x\varepsilon)\varepsilon$, $(x\varepsilon + y\varepsilon)^n = \sum_{k=0}^{n} \begin{bmatrix} n \\ k \end{bmatrix}_q (x\varepsilon)^k (y\varepsilon)^{n-k}$, $(x\varepsilon)^k = q^{\binom{k}{2}} x^k \varepsilon^k$ and

$$(x\varepsilon + y\varepsilon)^n = (x+y)\varepsilon (x+y)\varepsilon \cdots (x+y)\varepsilon = (x+y)(qx+y)(q^2x+y)\cdots(q^{n-1}x+y)\varepsilon^n.$$

**Remark**
For later use observe that comparing coefficients in

$$(y+qx)(y+q^2x)\cdots(y+q^n x) = \sum_{k=0}^{n} \begin{bmatrix} n \\ k \end{bmatrix}_q q^{\binom{k+1}{2}} x^k y^{n-k}$$

gives

$$\sum_{1 \le i_1 < i_2 < \cdots < i_k \le n} q^{i_1+i_2+\cdots+i_k} = q^{\binom{k+1}{2}} \begin{bmatrix} n \\ k \end{bmatrix}_q. \tag{50}$$

We shall also need the formula

$$\frac{1}{(1-x)(1-qx)\cdots(1-q^k x)} = \sum_{n \ge 0} \begin{bmatrix} n+k \\ n \end{bmatrix}_q x^n. \tag{51}$$

Let $G_k(x) = \sum_{n \geq 0} \begin{bmatrix} n+k \\ n \end{bmatrix}_q x^n$. Since $G_0(x) = \dfrac{1}{1-x}$ (51) follows from

$$(1-q^k x)G_k(x) = \sum_{n \geq 0} \begin{bmatrix} n+k \\ n \end{bmatrix}_q x^n - \sum_{n \geq 0} q^k \begin{bmatrix} n+k \\ n \end{bmatrix}_q x^{n+1} = \sum_{n \geq 0} \left( \begin{bmatrix} n+k \\ n \end{bmatrix}_q - q^k \begin{bmatrix} n+k-1 \\ n-1 \end{bmatrix}_q \right) x^n$$

$$= \sum_{n \geq 0} \left( \begin{bmatrix} n+k \\ n \end{bmatrix}_q - q^k \begin{bmatrix} n+k-1 \\ k \end{bmatrix}_q \right) x^n = \sum_{n \geq 0} \begin{bmatrix} n+k-1 \\ k-1 \end{bmatrix}_q x^n = G_{k-1}(x).$$

## 3. q-Fibonacci and q-Lucas polynomials

We will consider two different $q-$analogs of the Fibonacci polynomials.

### 3.1. Carlitz q-Fibonacci polynomials

These polynomials were introduced by L. Carlitz [5]. The $q-$analogs of the Fibonacci numbers $F_n(1,1,q)$ and $F_n(1,q,q)$ occurred earlier in I.Schur [15] in his polynomial versions of the Rogers-Ramanujan identities.

The arithmetical version can recursively be defined by

$$F_n(x,s,q) = xF_{n-1}(x,s,q) + q^{n-2}sF_{n-2}(x,s,q) \tag{52}$$
$$\text{with } F_0(x,s,q) = 0 \text{ and } F_1(x,s,q) = 1$$

and the combinatorial version by

$$f_n(x,s,q) = xf_{n-1}(x,s,q) + q^{n-1}sf_{n-2}(x,s,q) \tag{53}$$
$$\text{with } f_0(x,s,q) = 1 \text{ and } f_1(x,s,q) = x.$$

It is easily verified that

$$f_n(x,s,q) = \sum_{k=0}^{\lfloor \frac{n}{2} \rfloor} q^{k^2} \begin{bmatrix} n-k \\ k \end{bmatrix}_q s^k x^{n-2k}. \tag{54}$$

**Remark**

Equation (53) shows that the sequence $\left(f_n(x,s,q)\right)_{n \geq 0}$ is orthogonal.

(40) implies that for $s = -1$ the moments $c(2n, 0) = c_n(q)$ are $q-$analogs of the Catalan numbers $C_n$.

Computations show that

$$\left(c_n(q)\right)_{n \geq 0} = \left(1, 1, 1+q, 1+2q+q^2+q^3, 1+3q+3q^2+3q^3+2q^4+q^5+q^6, \cdots\right).$$

The $c_n(q)$ satisfy $c_n(q) = \sum_{k=0}^{n-1} q^k c_k(q) c_{n-1-k}(q)$ with $c_0(q) = 1$, but there is no closed formula.

Let us now give a combinatorial interpretation of $f_n(x,s,q)$.

For a Morse sequence $c$ of length $n$ which covers an interval $[m, m+n]$ with dots and dashes we define for $i \in \{0, \cdots, n-1\}$ a weight $v(i)$ by $v(i) = x$ if the interval $[i, i+1]$ is covered by a dot, $v(i) = q^i s$ if $[i, i+1]$ is the second part of a dash and $v(i) = 1$ else. Then

$$v(c) = \prod_{i=m}^{m+n-1} v(i) \tag{55}$$

and

$$f_n(x,s,q) = \sum_{c \in M_n} v(c). \tag{56}$$

The recurrence (53) results by considering the last element of $c$.
By considering the first element of $c$ we get another recurrence

$$f_n(x,s,q) = xf_{n-1}(x,qs,q) + qsf_{n-2}(x,q^2s,q) \tag{57}$$
with $f_0(x,s,q) = 1$ and $f_1(x,s,q) = x$.

Recurrence (52) for the arithmetical version $F_n(x,s,q)$ can be written in matrix form as

$$\begin{pmatrix} sF_{n-1}(x,qs,q) & F_n(x,s,q) \\ sF_n(x,qs,q) & F_{n+1}(x,s,q) \end{pmatrix} = \begin{pmatrix} 0 & 1 \\ q^{n-1}s & x \end{pmatrix} \begin{pmatrix} sF_{n-2}(x,qs,q) & F_{n-1}(x,s,q) \\ sF_{n-1}(x,qs,q) & F_n(x,s,q) \end{pmatrix} \tag{58}$$

Setting $U(x,s,q) = \begin{pmatrix} 0 & 1 \\ s & x \end{pmatrix}$ and $M_n(x,s,q) = \begin{pmatrix} sF_{n-1}(x,qs,q) & F_n(x,s,q) \\ sF_n(x,qs,q) & F_{n+1}(x,s,q) \end{pmatrix}$ we get

$$M_n(x,s,q) = U(x,q^{n-1}s,q)M_{n-1}(x,s,q) = U(x,q^{n-1}s,q)\cdots U(x,qs,q)(U(x,s,q). \tag{59}$$

Taking determinants, we get a $q$ – analog of **Cassini's identity**

$$F_n(x,s,q)F_n(x,qs,q) - F_{n-1}(x,qs.q)F_{n+1}(x,s,q) = q^{\binom{n}{2}}(-s)^{n-1}. \tag{60}$$

From
$$M_{m+n}(x,s,q) = M_m(x,q^ns,q)M_n(x,s,q)$$
$$= \begin{pmatrix} q^n sF_{m-1}(x,q^{n+1}s,q) & F_m(x,q^ns,q) \\ q^n sF_m(x,q^{n+1}s,q) & F_{m+1}(x,q^ns,q) \end{pmatrix} \begin{pmatrix} sF_{n-1}(x,qs,q) & F_n(x,s,q) \\ sF_n(x,qs,q) & F_{n+1}(x,s,q) \end{pmatrix}$$

we get

$$F_{m+n}(x,s,q) = F_m(x,q^ns,q)F_{n+1}(x,s,q) + q^n sF_{m-1}(x,q^{n+1}s,q)F_n(x,s,q). \tag{61}$$

There is also an analog of the doubling formula

$$f_{2n}(x,s,q) = \sum_{k=0}^{n} \begin{bmatrix} n \\ k \end{bmatrix}_q q^{kn}s^k x^{n-k} f_{n-k}(x,s,q). \tag{62}$$

Each Morse sequence $c = c_1 \cdots c_m$ of length $2n$ can be split into
$c = (c_1 \cdots c_{m-n})(c_{m-n+1} \cdots c_m) = de.$
If $e$ has $k$ dashes then its length is $n+k$ and therefore the length of $d$ is $n-k$.
Therefore, the sum over all $v(d)$ is $f_{n-k}(x,s,q)$.

It remains to show that the sum over all $v(e)$ is $\begin{bmatrix} n \\ k \end{bmatrix}_q q^{kn} s^k x^{n-k}$.

Since $v(a)v(b) = qv(b)v(a)$ this follows from the binomial theorem for $q-$commuting operators.

$$(v(a)+v(b))^n = \sum_{k=0}^n \begin{bmatrix} n \\ k \end{bmatrix}_q v(b)^k v(a)^{n-k} = \sum_{k=0}^n \begin{bmatrix} n \\ k \end{bmatrix}_q v(b\cdots b)x^{n-k} = \sum_{k=0}^n \begin{bmatrix} n \\ k \end{bmatrix}_q q^{k(n-k)+1+3+\cdots+2k-1}s^k x^{n-k}$$

$$= \sum_{k=0}^n \begin{bmatrix} n \\ k \end{bmatrix}_q q^{kn} s^k x^{n-k}.$$

For the generating function we get

$$\sum_{n\geq 0} f_n(x,s,q) z^n = \sum_{\ell \geq 0} z^\ell \sum_{k=0}^{\lfloor \ell/2 \rfloor} q^{k^2} \begin{bmatrix} \ell-k \\ k \end{bmatrix}_q s^k x^{\ell-2k} = \sum_{\substack{k,\ell \geq 0 \\ 2k \leq \ell}} q^{k^2} s^k z^{2k} \begin{bmatrix} \ell-k \\ k \end{bmatrix}_q x^{\ell-2k} z^{\ell-2k} \quad (63)$$

$$= \sum_{k\geq 0} q^{k^2} s^k z^{2k} \sum_{n\geq 0} \begin{bmatrix} n+k \\ k \end{bmatrix}_q x^n z^n = \sum_{k\geq 0} \frac{q^{k^2} s^k z^{2k}}{(1-xz)(1-qxz)\cdots(1-q^k xz)}$$

We can also define $q-$Lucas polynomials by

$$L_n(x,s,q) = tr M_n(x,s,q) = F_{n+1}(x,s,q) + s F_{n-1}(x,qs,q) = \sum_{k=0}^{\lfloor n/2 \rfloor} q^{k^2-k} \begin{bmatrix} n-k \\ k \end{bmatrix}_q \frac{[n]_q}{[n-k]_q} s^k x^{n-2k}.$$

(64)

**Remark.**
There are no nice recurrence relations for these $q-$Lucas polynomials. To remedy this failure H. Belbachir and A. Benmezai [3] proposed two alternatives. Let us call them $LB_n(x,s,q)$ and $\mathbf{LB}_n(x,s,q)$, by setting

$$\mathbf{LB}_n(x,s,q) = \sum_{k=0}^{\lfloor n/2 \rfloor} q^{k^2-k} \begin{bmatrix} n-k \\ k \end{bmatrix}_q \left(1 + \frac{[k]_q}{[n-k]_q}\right) s^k x^{n-2k},$$

$$LB_n(x,s,q) = \sum_{k=0}^{\lfloor n/2 \rfloor} q^{k^2} \begin{bmatrix} n-k \\ k \end{bmatrix}_q \left(1 + q^{n-2k} \frac{[k]_q}{[n-k]_q}\right) s^k x^{n-2k}.$$

(65)

For $q \to 1$ both converge to $L_n(x,s)$ and $LB_n(x,s,q) + \mathbf{LB}_n(x,s,q) = 2L_n(x,s,q)$.
The polynomials $\mathbf{LB}_n(x,s,q)$, satisfy

$$\mathbf{LB}_n(x,s,q) = x\mathbf{LB}_{n-1}(x,s,q) + q^{n-2} s \mathbf{LB}_{n-2}(x,s,q), \quad (66)$$

and the polynomials $LB_n(x,s,q)$ satisfy

$$LB_n(x,s,q) = x LB_{n-1}(x,qs,q) + qs LB_{n-2}(x,q^2 s,q). \quad (67)$$

These are nice companions to (53) and (57).

### 3.2. The polynomials $Fib_n(x,s,q)$ and $Luc_n(x,s,q)$.

For a Morse sequence $c = c_1 c_2 \cdots c_m$ with dashes at $i_1, i_2, \cdots, i_k$
we define the weight
$$W(c,x,s) = q^{i_1+i_2+\cdots+i_k} s^k x^{m-k} \tag{68}$$
and define the arithmetic $q$ – Fibonacci polynomial $Fib_n(x,s,q)$ by
$$Fib_n(x,s,q) = \sum_{c \in M_{n-1}} W(c,x,s). \tag{69}$$
Since $W(c_1 c_2 \cdots c_m, x, s) = W(c_1, x, s) W(c_2 \cdots c_m, x, qs)$ we get
$$Fib_n(x,s,q) = xFib_{n-1}(x,qs,q) + qsFib_{n-2}(x,qs,q) \tag{70}$$
with initial values $Fib_0(x,s,q) = 0$ and $Fib_1(x,s,q) := 1$.
Since a Morse sequence of length $n-1$ with $k$ dashes has $n-1-k$ elements formula (50) gives
$$Fib_n(x,s,q) = \sum_{k=0}^{\lfloor \frac{n-1}{2} \rfloor} \begin{bmatrix} n-1-k \\ k \end{bmatrix}_q q^{\binom{k+1}{2}} s^k x^{n-1-2k}. \tag{71}$$
The first terms are $0, 1, x, x^2+qs, x^3+(1+q)qsx, x^4+(1+q+q^2)qsx^2+q^3s^2$.

If we consider the last element of a Morse sequence we get the recurrence
$$Fib_n(x,s,q) = xFib_{n-1}(x,s,q) + q^{n-2}sFib_{n-2}\left(x,\frac{s}{q},q\right). \tag{72}$$
If the last element is a dot we get the first term of the recurrence. For the second term we consider all Morse sequences with length $n-1$ and $k$ dashes. They have $n-k-1$ elements. The weight of those sequences where the last element is a dash is by (50)
$$q^{n-1-k} s \begin{bmatrix} n-3-k+1 \\ k-1 \end{bmatrix}_q q^{\binom{k}{2}} s^{k-1} x^{n-1-2k} = q^{n-2} s \begin{bmatrix} n-2-k \\ k-1 \end{bmatrix}_q q^{\binom{k}{2}} s^{k-1} x^{n-1-2k}.$$
Summing over all $k$ gives the second term of (72).
The polynomials $Fib_n(x,s,q)$ are not orthogonal but satisfy the recurrence (cf. [9])

$$Fib_n(x,s,q) = xFib_{n-1}(x,s,q) + q^{n-2}sxFib_{n-3}(x,s,q) + q^{n-2}s^2 Fib_{n-4}(x,s,q), \tag{73}$$

because by (72) we get $Fib_n(x,s,q) = xFib_{n-1}(x,s,q) + q^{n-2}sFib_{n-2}\left(x,\frac{s}{q},q\right)$ and by replacing

$n \to n-2$ and $s \to \frac{s}{q}$ in (70) $Fib_{n-2}\left(x,\frac{s}{q},q\right) = xFib_{n-3}(x,s,q) + sFib_{n-4}(x,s,q)$.

As special case we get an analog of (23) (cf. [9]):

$$\left(Fib_n\left(1,-\frac{1}{q},q\right)\right)_{n \geq 0} = \left(0,1,1,0,-q,-q^2,0,q^5,q^7,0,-q^{12},-q^{15},\cdots\right). \tag{74}$$

Denoting by $r(n) = \dfrac{n(3n-1)}{2}$ the pentagonal numbers this means

$$Fib_{3n}\left(1, -\frac{1}{q}\right) = 0,$$

$$Fib_{3n+1}\left(1, -\frac{1}{q}\right) = (-1)^n q^{r(n)}, \tag{75}$$

$$Fib_{3n+2}\left(1, -\frac{1}{q}\right) = (-1)^n q^{r(-n)}.$$

**Proof.**

For $(x,s) = \left(1, -\dfrac{1}{q}\right)$ (73) reduces to

$$Fib_n\left(1, -\frac{1}{q}, q\right) = Fib_{n-1}\left(1, -\frac{1}{q}, q\right) - q^{n-3} Fib_{n-3}\left(1, -\frac{1}{q}, q\right) + q^{n-4} Fib_{n-4}\left(1, -\frac{1}{q}, q\right).$$

Therefore $a_n = Fib_n\left(1, -\dfrac{1}{q}, q\right) + q^{n-3} Fib_{n-3}\left(1, -\dfrac{1}{q}, q\right)$ satisfies $a_n = a_{n-1}$.

Since $Fib_3(x, s, q) = x^2 + qs$ we get $a_3 = 0$ and thus also $a_n = 0$ for $n \geq 3$.

This gives

$$Fib_n\left(1, -\frac{1}{q}, q\right) = -q^{n-3} Fib_{n-3}\left(1, -\frac{1}{q}, q\right). \tag{76}$$

(75) is true for $n = 0$. Identities (75) follow by induction from $r(n) - r(n-1) = 3n - 2$,

$r(-n) - r(-n+1) = 3n - 1$, $Fib_{3n+1}\left(1, -\dfrac{1}{q}, q\right) = -q^{3n-2} Fib_{3n-2}\left(1, -\dfrac{1}{q}, q\right)$ and

$$Fib_{3n+2}\left(1, -\frac{1}{q}, q\right) = -q^{3n-1} Fib_{3n-1}\left(1, -\frac{1}{q}, q\right).$$

There is also a recurrence of a rather curious kind. If we denote by $D_q$ the $q$-differential operator

$$D_q p(x) = \frac{p(x) - p(qx)}{(1-q)x} \tag{77}$$

then we get

$$Fib_n(x, s, q) = x Fib_{n-1}(x, s, q) + (q-1)s D_q Fib_{n-1}(x, s, q) + s Fib_{n-2}(x, s, q). \tag{78}$$

Using

$$(q-1)D_q Fib_n(x,s,q) = \sum_{k=0}^{\lfloor \frac{n-1}{2} \rfloor} \begin{bmatrix} n-1-k \\ k \end{bmatrix}_q q^{\binom{k+1}{2}} s^k \left(q^{n-1-2k} - 1\right) x^{n-2-2k}$$

$$= \sum_{k=0}^{\lfloor \frac{n-1}{2} \rfloor} \begin{bmatrix} n-2-k \\ k \end{bmatrix}_q q^{\binom{k+1}{2}} s^k \left(q^{n-1-k} - 1\right) x^{n-2-2k} = q^{n-1} Fib_{n-1}\left(x, \frac{s}{q}, q\right) - Fib_{n-1}(x,s,q)$$

and (72) gives (78).

Formula (78) implies that

$$Fib_n(x,s,q) = F_n\left(x + (q-1)sD_q, s\right)1. \tag{79}$$

To see this, define the linear operator $T_s$ on the polynomials $p(x)$ by

$$T_s x^n = \left(x + (q-1)sD_q\right)^n 1. \tag{80}$$

Applying $T_s$ to the identity $F_n(x,s) = xF_{n-1}(x,s) + sF_{n-2}(x,s)$ we get
$T_s F_n(x,s) = (x + (q-1)sD_q)T_s F_{n-1}(x,s) + sT_s F_{n-2}(x,s)$. Since $T_s(0) = 0$ and $T_s(1) = 1$ we get
$T_s F_n(x,s) = Fib_n(x,s,q)$.

We define the corresponding $q$ – **Lucas polynomials** by

$$Luc_n(x,s,q) = L_n\left(x + (q-1)sD_q, s\right)1. \tag{81}$$

The first terms are

$2, \; x, \; x^2 + (1+q)s, \; x^3 + \left(1+q+q^2\right)sx, \; x^4 + \left(1+q+q^2+q^3\right)sx^2 + q\left(1+q^2\right)s^2.$

From $L_n(x,s) = F_{n+1}(x,s) + sF_{n-1}(x,s)$ we get by applying the linear operator $T_s$
the identity

$$Luc_n(x,s,q) = Fib_{n+1}(x,s,q) + sFib_{n-1}(x,s,q). \tag{82}$$

This implies the explicit formula

$$Luc_n(x,s,q) = \sum_{j=0}^{\lfloor \frac{n}{2} \rfloor} q^{\binom{j}{2}} \begin{bmatrix} n-j \\ j \end{bmatrix}_q \frac{[n]_q}{[n-j]_q} s^j x^{n-2j}. \tag{83}$$

For the proof compare the coefficient of $s^j x^{n-2j}$. The right-hand side gives

$$q^{\binom{j+1}{2}} \begin{bmatrix} n-j \\ j \end{bmatrix}_q + q^{\binom{j}{2}} \begin{bmatrix} n-1-j \\ j-1 \end{bmatrix}_q = q^{\binom{j}{2}} \begin{bmatrix} n-j \\ j \end{bmatrix}_q \left(q^j + \frac{[j]_q}{[n-j]_q}\right) = q^{\binom{j}{2}} \begin{bmatrix} n-j \\ j \end{bmatrix}_q \frac{[n]_q}{[n-j]_2}.$$

The combinatorial analog satisfies

$$\begin{aligned} luc_n(x,s,q) &= Luc_n(x,s,q) \text{ for } n > 0 \\ luc_0(x,s,q) &= 1. \end{aligned} \tag{84}$$

Thus all $luc_n(x,s,q)$ are monic polynomials with degree $\deg_x luc_n(x,s) = n$.

As analog of (38) we get ([4])

$$\sum_{k=0}^{\lfloor \frac{n}{2} \rfloor} (-s)^k \left( \begin{bmatrix} n \\ k \end{bmatrix}_q - \begin{bmatrix} n \\ k-1 \end{bmatrix}_q \right) fib_{n-2k}(x,s,q) = x^n,$$

$$\sum_{k=0}^{\lfloor \frac{n}{2} \rfloor} (-s)^k \begin{bmatrix} n \\ k \end{bmatrix}_q luc_{n-2k}(x,s,q) = x^n.$$

(85)

As a corollary we see that

$$\Lambda_{Fib(x,s,q)}(x^{2n}) = (-qs)^n \frac{1}{[n+1]_q} \begin{bmatrix} 2n \\ n \end{bmatrix}_q \tag{86}$$

and

$$\Lambda_{Luc(x,s,q)}(x^{2n}) = (-s)^n \begin{bmatrix} 2n \\ n \end{bmatrix}_q. \tag{87}$$

**Proof**

The first identity follows from the second by using (82).
For the proof of the second one we consider (cf. [12])

$$H_n(x,s,q) = \sum_{k=0}^{n} (-s)^k \begin{bmatrix} n \\ k \end{bmatrix}_q l_{n-2k}(x,s). \tag{88}$$

We get

$$H_n(x,s,q) = xH_{n-1}(x,s,q) + s(1-q^{n-1})H_{n-2}(x,s,q) \tag{89}$$

and

$$H_n\big(x+(q-1)sD_q, s, q\big)1 = x^n, \tag{90}$$

which implies (85).

Let us recall that $l_n(x,s) = xl_{n-1}(x,s) + t_n(s)l_{n-2}(x,s)$ with $t_n(s) = s$ for $n > 0$ and $t_0(s) = 2s$.

For odd $n$ (89) follows from

$$H_n(x,s,q) = \sum_{k=0}^{\lfloor \frac{n}{2} \rfloor} (-s)^k \left( \begin{bmatrix} n-1 \\ k \end{bmatrix}_q + q^{n-k} \begin{bmatrix} n-1 \\ k-1 \end{bmatrix}_q \right) l_{n-2k}(x,s)$$

$$= \sum_k (-s)^k \begin{bmatrix} n-1 \\ k \end{bmatrix}_q (l_{n-2k}(x,s) - sl_{n-2k-2}(x,s)) + s\sum_k l_{n-2k-2}(x,s)(-s)^k (1-q^{n-k-1}) \begin{bmatrix} n-1 \\ k \end{bmatrix}_q$$

$$= \sum_k (-s)^k \begin{bmatrix} n-1 \\ k \end{bmatrix}_q xl_{n-1-2k}(x,s) + s\sum_k l_{n-2k-2}(x,s)(-s)^k (1-q^{n-1}) \begin{bmatrix} n-2 \\ k \end{bmatrix}_q$$

$$= xH_{n-1}(x,s,q) + s(1-q^{n-1})H_{n-2}(x,s,q).$$

For even $n = 2m$ we get

$$H_n(x,s,q) = \sum_{k=0}^{m}(-s)^k\left(\begin{bmatrix}n-1\\k\end{bmatrix}_q + q^{n-k}\begin{bmatrix}n-1\\k-1\end{bmatrix}_q\right)l_{n-2k}(x,s)$$

$$= (-s)^m\begin{bmatrix}n-1\\m\end{bmatrix}_q + \sum_{k=0}^{m-1}(-s)^k\begin{bmatrix}n-1\\k\end{bmatrix}_q \left(l_{n-2k}(x,s) - sl_{n-2k-2}(x,s)\right) + s\sum_{k=0}^{m-1}l_{n-2k-2}(x,s)(-s)^k\left(1-q^{n-k-1}\right)\begin{bmatrix}n-1\\k\end{bmatrix}_q$$

$$= (-s)^m\begin{bmatrix}n-1\\m\end{bmatrix}_q + \sum_{k=0}^{m-1}(-s)^k\begin{bmatrix}n-1\\k\end{bmatrix}_q \left(l_{n-2k}(x,s) - t_{n-2k-2}(s)l_{n-2k-2}(x,s)\right) - (-s)^m\begin{bmatrix}n-1\\m-1\end{bmatrix}_q$$

$$+ s\sum_{k=0}^{m-1}l_{n-2k-2}(x,s)(-s)^k\left(1-q^{n-1}\right)\begin{bmatrix}n-2\\k\end{bmatrix}_q$$

$$= \sum_{k=0}^{m-1}(-s)^k\begin{bmatrix}n-1\\k\end{bmatrix}_q xl_{n-2k}(x,s) + s\left(1-q^{n-1}\right)\sum_{k=0}^{m-1}l_{n-2k-2}(x,s)(-s)^k\begin{bmatrix}n-2\\k\end{bmatrix}_q$$

$$= xH_{n-1}(x,s,q) + s\left(1-q^{n-1}\right)H_{n-2}(x,s,q).$$

Identity (90) is trivially true for $n=0$ and $n=1$.
By induction we get
$$H_n(x+(q-1)sD_q, s, q)1 = \left(x+(q-1)sD_q\right)x^{n-1} + s\left(1-q^{n-1}\right)x^{n-2} = x^n.$$

It is interesting to note that (47) implies
$$H_n(x,s,q) = \sum_{k=0}^{n}\begin{bmatrix}n\\k\end{bmatrix}_q \alpha(x,s)^k \beta(x,s)^{n-k} = R_n(\alpha(x,s), \beta(x,s), q). \tag{91}$$

For the generating function we get

$$\sum_{n\geq 0} fib_n(x,s,q)z^n = \sum_{k\geq 0} \frac{q^{\binom{k+1}{2}}s^k z^{2k}}{(1-xz)(1-qxz)\cdots(1-q^k xz)}. \tag{92}$$

**Proof.**

$$\sum_{n\geq 0} fib_n(x,s,q)z^n = \sum_{\ell \geq 0} z^\ell \sum_{k=0}^{\lfloor \ell/2 \rfloor} q^{\binom{k+1}{2}}\begin{bmatrix}\ell-k\\k\end{bmatrix}_q s^k x^{\ell-2k} = \sum_{\substack{k,\ell \geq 0 \\ 2k \leq \ell}} q^{\binom{k+1}{2}} s^k z^{2k}\begin{bmatrix}\ell-k\\k\end{bmatrix}_q x^{\ell-2k} z^{\ell-2k}$$

$$= \sum_{k\geq 0} q^{\binom{k+1}{2}} s^k z^{2k} \sum_{n\geq 0}\begin{bmatrix}n+k\\k\end{bmatrix}_q x^n z^n = \sum_{k\geq 0} \frac{q^{\binom{k+1}{2}} s^k z^{2k}}{(1-xz)(1-qxz)\cdots(1-q^k xz)}.$$

**Remark**
There are nice $q-$analogs of (22) and (36). We state them without proof. Proofs can be found in [4] and [10].

Let $R_n(x,y,q) = \sum_{k=0}^{n}\begin{bmatrix}n\\k\end{bmatrix}_q x^k y^{n-k}$ be the Rogers-Szegö polynomials. Then

$$\sum_{k=0}^{\left\lfloor\frac{n-1}{2}\right\rfloor} q^{\binom{k+1}{2}} \begin{bmatrix} n-1-k \\ k \end{bmatrix}_q (-xy)^k R_{n-1-2k}(x,y) = \frac{x^n - y^n}{x - y},$$

$$\sum_{k=0}^{\left\lfloor\frac{n}{2}\right\rfloor} q^{\binom{k}{2}} \begin{bmatrix} n-k \\ k \end{bmatrix}_q \frac{[n]_q}{[n-k]_q} (-xy)^k R_{n-2k}(x,y) = x^n + y^n.$$
(93)

### 4. Concluding remarks

We have seen that the two classes of $q-$Fibonacci and q-Lucas polynomials differ in many aspects. But the formulas $f_n(x,s,q) = \sum_{k=0}^{\left\lfloor\frac{n}{2}\right\rfloor} q^{k^2} \begin{bmatrix} n-k \\ k \end{bmatrix}_q s^k x^{n-2k}$ and

$fib_n(x,s,q) = \sum_{k=0}^{\left\lfloor\frac{n}{2}\right\rfloor} \begin{bmatrix} n-k \\ k \end{bmatrix}_q q^{\binom{k+1}{2}} s^k x^{n-2k}$ show that the linear map $\Phi$ from the polynomials in $x,s$ into itself defined by

$$\Phi(x^n s^k) = q^{\binom{k}{2}} x^n s^k$$
(94)

satisfies

$$\Phi(fib_n(x,s,q)) = f_n(x,s,q).$$
(95)

Moreover, we get by comparing coefficients

$$\Phi(s^k fib_n(x,s,q)) = q^{\binom{k}{2}} s^k f_n(x, q^k s, q),$$
$$\Phi(s^k luc_n(x,s,q)) = q^{\binom{k}{2}} s^k l_n(x, q^k s, q).$$
(96)

Therefore some identities for one class can be transferred to the other class.

For example, the recurrence (52) translates into (72)
and (85) implies analogous results for the Carlitz polynomials

$$\sum_{k=0}^{\left\lfloor\frac{n}{2}\right\rfloor} (-s)^k q^{\binom{k}{2}} \left( \begin{bmatrix} n \\ k \end{bmatrix}_q - \begin{bmatrix} n \\ k-1 \end{bmatrix}_q \right) f_{n-2k}(x, q^k s, q) = x^n,$$

$$\sum_{k=0}^{\left\lfloor\frac{n}{2}\right\rfloor} (-s)^k q^{\binom{k}{2}} \begin{bmatrix} n \\ k \end{bmatrix}_q l_{n-2k}(x, q^k s, q) = x^n.$$
(97)